\documentclass[11pt]{article}
\usepackage{gensymb}

\usepackage{fullpage,amssymb,amsmath}
\usepackage{graphicx}
\usepackage{natbib}
\usepackage[capitalise]{cleveref}
\newcommand{\dd}{\mathrm{d}}

\title{The optimal way of folding a bond notebook page into a bookmark}
\author{Chenguang Zhang\thanks{Department of Mathematics, Massachusetts Institute of Technology. Email: \texttt{cgzhang@mit.edu}.}}

\begin{document}
\maketitle

\begin{abstract}
The article explores a simple yet interesting geometric problem: how to fold a notebook page so it is the most visible when the notebook is closed (i.e., as a bookmark). We start with a square page, then proceed to the rectangular page. When an additional constraint is added to limit the vertical extent of the rectangular page, the problem shows the first-order phase transition: the optimal way of folding changes drastically when the aspect ratio of the page exceeds the critical value of about $1.20711$.
\end{abstract}

\section{Introduction}\label{introduction}
\begin{figure}[h]
	\centering
	\includegraphics[width=0.6\linewidth]{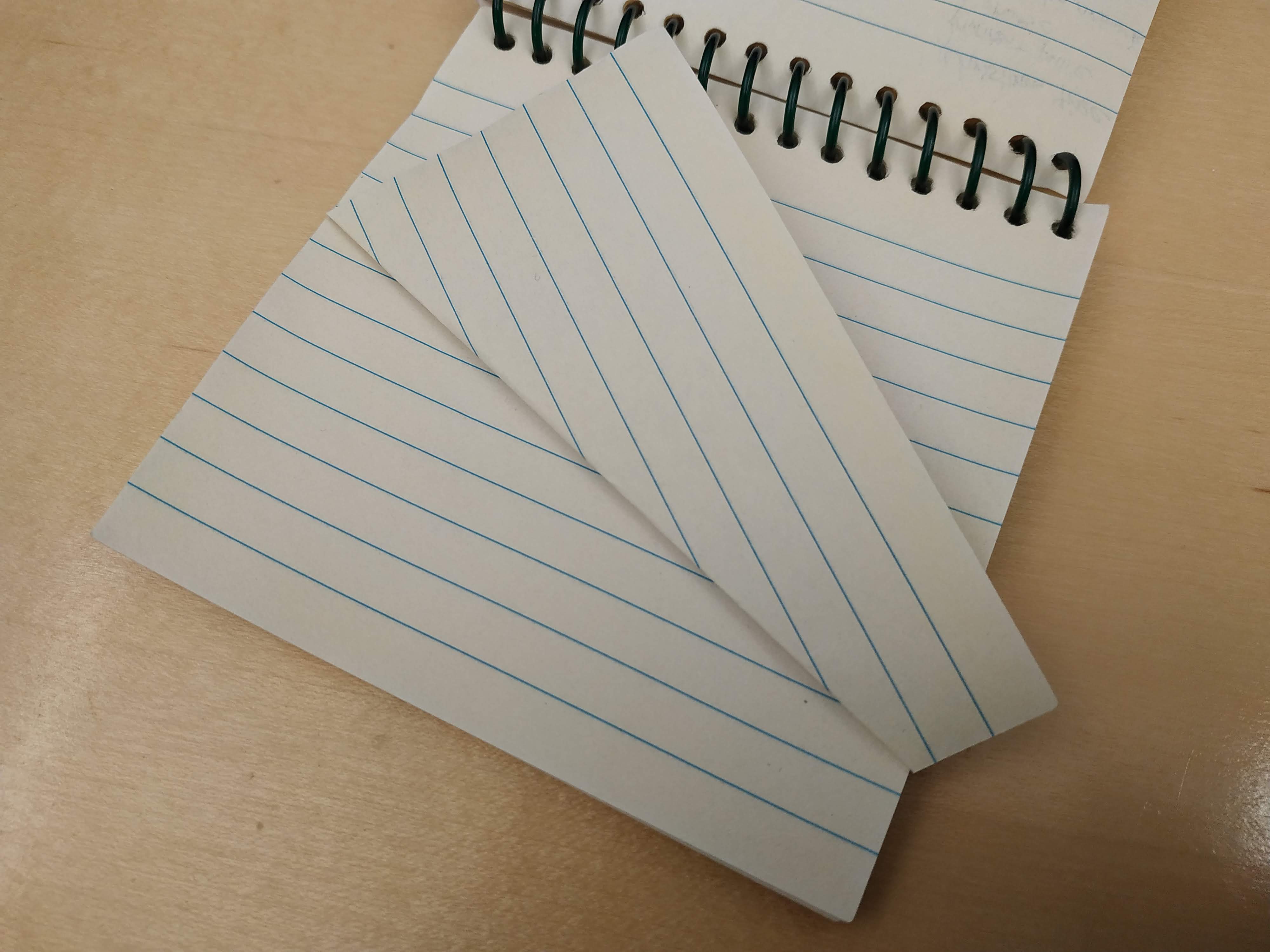}
	\caption{The square notebook that inspired this work. This fold here is actually a close approximation of the optimal solution (\cref{sec:case2}).}
	\label{fig:illu}
\end{figure}

This problem is inspired by a small square notebook I am using recently (\cref{fig:illu}). One day I decided to use a page as bookmark: by folding it sideways, I can easily locate it when the notebook is closed. \cref{fig:illu} shows the two requirements of a solution. First, the folded page must exceed the boundary of the notebook. Second, all pages are bond from the top and cannot be torn off. 

This problem is related with the rich field of origami, which I very briefly review here. Origami is the art of paper folding that has entertained mankind for a long time. Since the last century, it has seen significant developments that, interestingly, were far beyond its originally recreational and artistic nature. 

From the theoretical side, mathematicians have identified the complete set of  axioms of origami: the Huzita-Hatori axioms \citep{alperin2009,lang2003}.
The application of rigorous mathematical formulations also helped creating a great number of new origami designs \citep{lang2011}. Although innocent at the first glance, origami surpasses the power of ``compass and straighthedge'' and can solve third-order mathematical problems including the ``angle trisection'' and ``doubling the cube''. Readers can find ample examples of mathematical origami in \citet{hull2012}.

From the technological side, origami is a \emph{generic} methodology to transform between 2d and 3d geometries. In the forwared direction (2d$\to$3d), it inspires the designs of toys,
home appliances, and medical devices such as stent. In the backwared direction (3d$\to$2d), it is used in airbag and solar panels, or wherever a thin and large structure needs to be folded and compactly stored, yet remains easily deployable. A celebrated example is the Miura fold created by Koryo Miura (a Japanese astrophysicist), which can fold a solar panel to about 1/25 of its original area, and later deployed by a single pulling action. The traditional origami and many of its engineering incarnations use mechanical loads to fold/unfold, which can be impossible, for example, at nano-scale or within human blood vessel. For this reason, researchers developed the so-called active origami, where folding/unfolding are achieved by temperature, chemicals, etc. It is now an active research area \citep[see][chap. 1]{hernandez2018}.

Although related with origami, the problem considered in this paper has two key differences. The first is that one side of the page is constrained (i.e., fixed). The second is we only have a single fold (although nothing prevents multiple folds from being considered in the future). For these reasons, the author would rather call it the bookmark problem than the attractive yet misleading ``constrained origami''.

\section{The square page}\label{sec:analysis}
Referring to \cref{fig:illu}, it is clear that we have an optimization problem that can be more formally stated as: given the unit square $[0,1]\times[0,1]$, the top edge of which is fixed, find the optimal fold that maximizes $x_e$: $x_e$ being the maximum $x$-coordinate of the folded shape. 

Note that without any fold (the null fold), we have $x_e=1$. In this paper we assume, without loss of generality, that the fold is towards upper-right.
Mathematically, a fold operation is to first choose a straight line (the crease) passing through the square, then reflect part of the square about it. To facilitate description, we make the following definitions:
\begin{description}
	\item[extreme point] the point with $x = x_e$. It is the rightmost point on the folded shape.
	\item[folded part] the part of square gets reflected. It is the part that gets temporarily detached from your table plane during folding.
\end{description}

Before going into details, we can make one general observation of the problem: $x_e$ must occur at one corner of the folded part. The reason is that if $x_e$ occurs at some other point $\boldsymbol{p}$ (either on the edge or within the folded part), then $\boldsymbol{p}$ will have certain neighboring points with greater $x$; this leads to contradiction. 

\begin{figure}[h]
	\centering
	\includegraphics[width=0.75\linewidth]{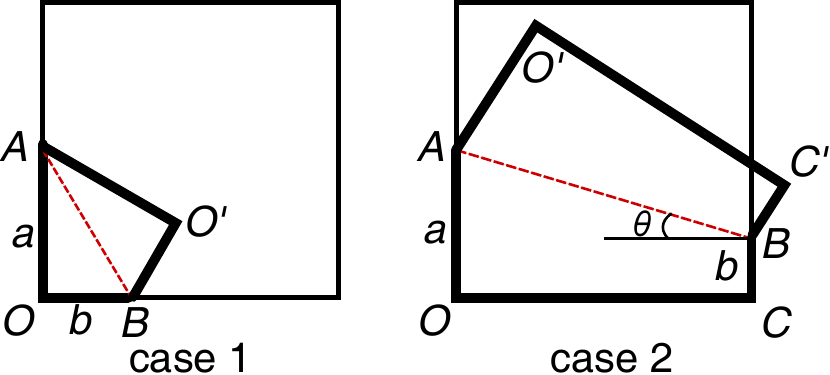}
	\caption{The two cases of folding a square towards upper-right when the top edge is fixed. The red dashed lines mark the creases. The edge length of the folded part is represented by lower case letters. The location of folded part after folding is drawn as well, with new corner locations notated by additional primes.}
	\label{fig:cases}
\end{figure}

Because the fold is towards upper-right, there are only two possible ways depending on how the crease crosses the square edges (\cref{fig:cases}). They are referred to as case 1 and 2 and discussed in turn. 

\subsection{Case 1}\label{sec:case1}
We first discuss when the crease passes through the left and bottom edges. In this case the folded part is $\triangle OAB$, so $x_e$ occurs either at $O'$ or $B$. Because the maximum possible $x_B$ is 1, $B$ does not quality a bookmark, we only focus on $O'$.

To find $x_{O'}$, notice that $\triangle OAB$ has height 
$$
h=\frac{ab}{\sqrt{a^2+b^2}}
$$
The length of $OO'$ is $2h$, a multiplication with $\cos (\angle BOO')=\cos (\angle OAB)$ gives
\begin{equation}\label{eqn:case1}
x_e(a,b)=x_{O'}(a,b) = \frac{2a^2b}{{a^2+b^2}}
\end{equation}

To find the maximum of \cref{eqn:case1}, we notice that $y/(y+\alpha)$ ($\alpha>0$) is monotonically increasing, thus we can directly set $a=1$ in \cref{eqn:case1}. The function becomes
\begin{equation}
x_{e}(1,b) = \frac{2b}{{1+b^2}}.
\end{equation}
It is clearly limited by the well-known inequality $\alpha^2+\beta^2\ge 2\alpha\beta$, and has a maximum of 1. Thus this fold never qualifies a bookmark.

\subsection{Case 2}\label{sec:case2}
Here the folded part is a trapzoid (\cref{fig:cases}). For the fold to point upper-right, we have $a\ge b$. Consider angle $\theta$ in the figure, $\tan\theta = a-b$. Because $0\le a-b\le 1$, $\theta \le 45\degree$. As a result, $O'$ never exceeds $C'$ (i.e., $O'C'$ never becomes vertical or tilts rightward). So we only need to focus on $x_{C'}$. For simplicity, we calculate the ``excess'' defined as $e = x_e-1$, which represents how much of the folded page remains visible after the notebook is closed. From basic trigonometry

$$e = 2 b \cos\theta \sin\theta = b\sin 2\theta$$

Using $\tan\theta = a-b$ and the tangent half-angle formula, we have
\begin{equation}\label{eqn:case2}
e(a,b) = \frac{2b(a-b)}{1+(a-b)^2}
\end{equation}

\begin{figure}[h]
	\centering
	\includegraphics[width=0.5\linewidth]{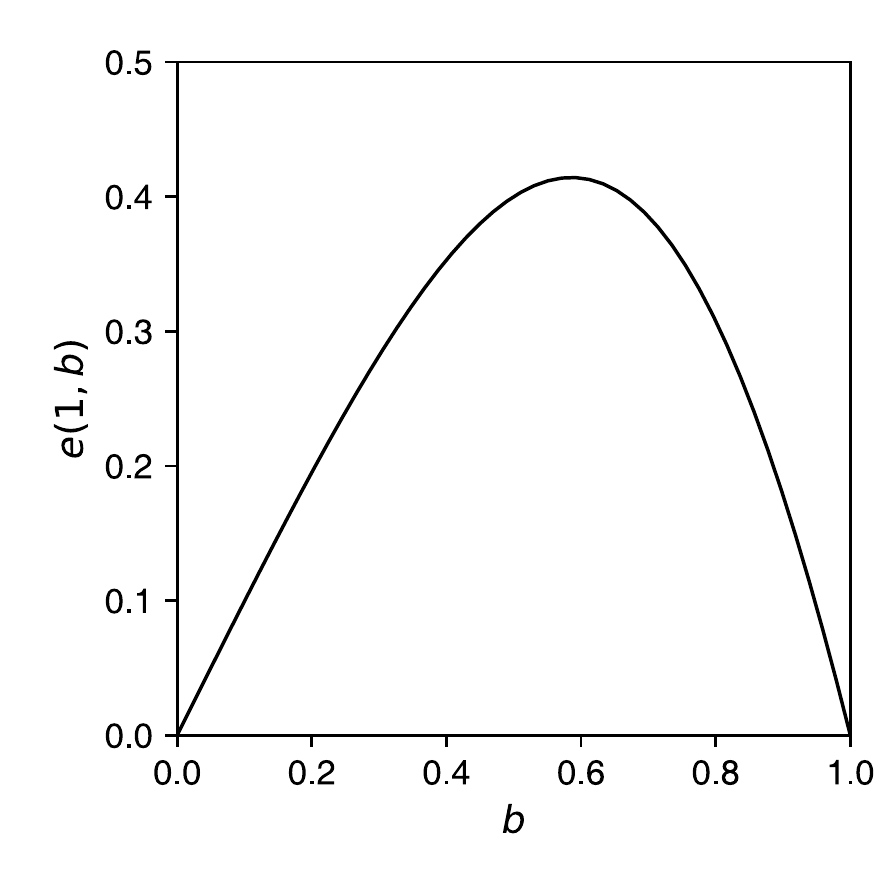}
	\caption{The curve of $e(b)$ in \cref{eqn:case2a}.}
	\label{fig:eb}
\end{figure}

This function is more complex, but we can use that $c/(1+c^2)$ is monotonically increasing for $0\le c\le 1$, and thus $\partial e/\partial a = \partial e/\partial c \cdot \partial c/\partial a>0$. So again we set $a=1$, the resultant function
\begin{equation}\label{eqn:case2a}
e(1,b) = \frac{2b(1-b)}{1+(1-b)^2}
\end{equation}
is no longer monotonic (\cref{fig:eb}). We first test it against some specific values:
\begin{enumerate}
	\item $b=0$ gives $e=0$. In fact, $e(a,0)=0$ regardless of $a$.
	\item $b=1$ gives $e=0$, because this fold just flips the page upwards. In fact, $e=0$ whenever $b=a$ (i.e., fold upwards).
\end{enumerate}
The optimum of \cref{eqn:case2a} is achieved when $\dd e(1,b)/\dd b = 0$, which gives $b = 2-\sqrt{2}\approx0.586$, $e=\sqrt{2}-1\approx0.414$. This means to have the most visible bookmark: first pick the top-left corner, then pick from the right edge a point $58.6\%$ from the bottom, then fold.

\subsection{Fold not exceeding top edge}\label{sec:foldconstraint}
One practical issue with the above solution is that the folded page exceeds its top edge. Clearly the maximum $y$, or $y_e$, occurs at $O'$:
\begin{equation}\label{eqn:ye}
y_e(a,b)=y_{O'}(a,b) = \frac{2a}{1+(a-b)^2}.
\end{equation}
Given the optimal solution $a =1, b= 2-\sqrt{2}$ of the previous section, we have
$$
y_e = 1 + \frac{\sqrt{2}}{2} \approx 1.707
$$
This means when the notebook is closed (after the page is fold), the page is passively folded again. For those who find this undesirable, we should limit the maximum $y$ with an additional constraint $y_e\le1$ and find a new solution. Now the problem can be succinctly presented as
\begin{equation}\label{eqn:case2y}
\begin{array}{rrclcl}
\displaystyle \max_{a,b} & \multicolumn{3}{l}{   \displaystyle e(a,b)=\frac{2b(a-b)}{1+(a-b)^2}  } \\
\textrm{s.t.} & b & \le & a \\
              & 2a & \le & 1+(a-b)^2 \\
              & a,b & \ge & 0
\end{array}
\end{equation}

\begin{figure}[h]
	\centering
	\includegraphics[width=0.5\linewidth]{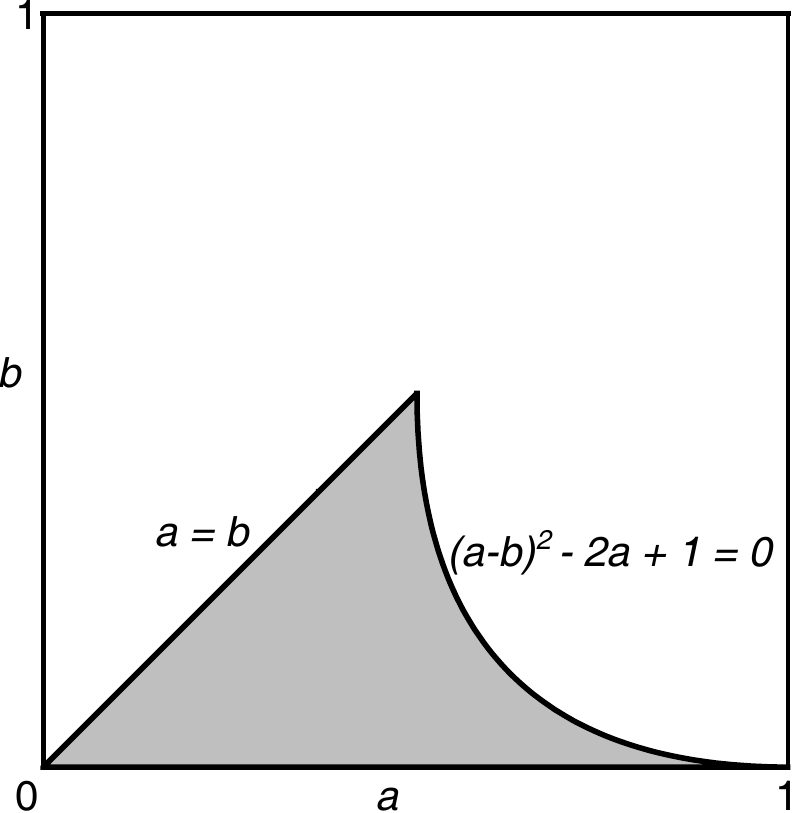}
	\caption{The feasible region defined by the constraints in \cref{eqn:case2y}.}
	\label{fig:region}
\end{figure}

It is important to realize that we are maximizing the $x$ of point $C'$ subject to the constraint on the $y$ of point $O'$, apparently due to that the two points are linked together. \cref{fig:region} shows the ``feasible region'' defined by the two constraints. One of the boundary is nonlinear, which makes the problem more complex. Still, the fact that $\partial e/\partial a\ge 0$ means that the maximum of $e$ occurs at the nonlinear boundary ($e=0$ along the other two straight boundaries). From the nonlinear boundary equation, we obtain
$a=1+b-\sqrt{2b}$ (the second solution $1+b+\sqrt{2b}$ is discarded), substitute into the objective function gives
$$
\max_{b} \, \left(b-\frac{b^2}{1+b-\sqrt{2b}}\right), b\in\left[0,\frac{\sqrt{2}}{2}\right]
$$
The upper limit of $b$ occurs at the intersection of the nonlinear boundary and the boundary $a=b$.
Apply $\dd /\dd b=0$ gives a 3rd order algebraic equation in terms of $\sqrt{b}$
$$\sqrt{b}^3 - 2\sqrt{2}\sqrt{b}^2+4\sqrt{b}-\sqrt{2}=0.$$
The only real solution of which is
$$
b = \frac{\sqrt[3]{7 \sqrt{33}+9}}{\sqrt[3]{9}}-\frac{8}{\sqrt[3]{3 \left(7 \sqrt{33}+9\right)}}\approx 0.248,
$$
which in turn gives $a \approx 0.543$. Note that here we must round to smaller numerical values to avoid violating the constraint $y_e\le1$. All together we have the optimal fold
$$e(0.543, 0.248) = 0.135.$$
The maximum $y$ value is $y_e = 0.999057$. As expected, it is close to but less than 1. By the way, it is easier to directly use numerical optimization to find the optimal solution (e.g, \texttt{FindMaximum} in Mathematca). 

Before ending this section, I admit that this solution is not easy to pull off in reality. The rough approximation $e(0.5, 0.25)=0.118$ is 12.6\% less than the optimum; not bad!
\subsection{A summary picture}
\begin{figure}[h]
	\centering
	\includegraphics[width=0.6\linewidth]{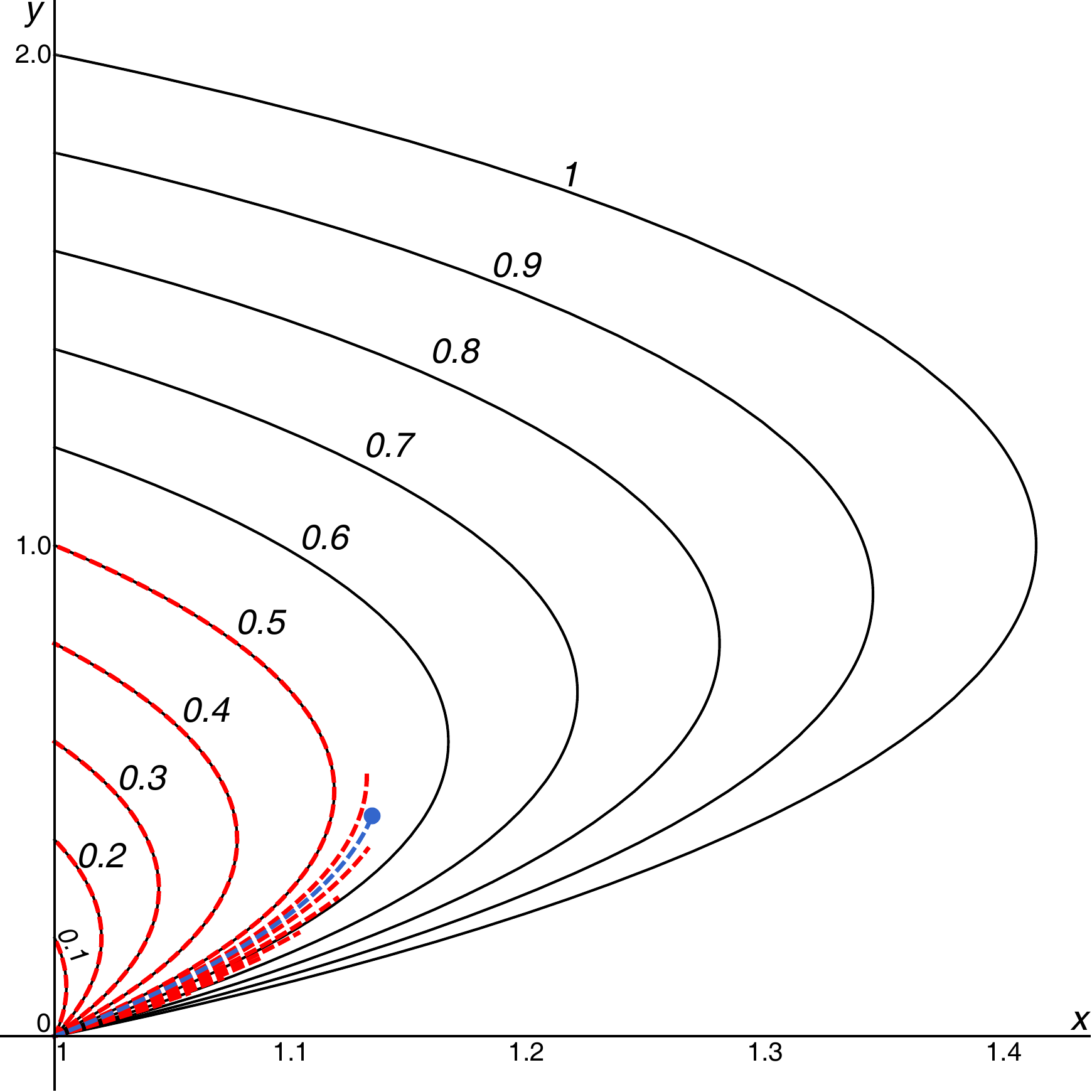}
	\caption{The trajectories of the bottom-right corner of the square in case 2  when we fix $a$ but vary $b$. Most curves are labeled by their $a$ values. The black curves are from the unconstrained folds (the top edge can be exceeded), the red curves are from the constrained folds (the top edge cannot be exceeded). Further, the curve with $a=0.543$ is shown in blue, with the optimal solution represented by a blue dot.}
	\label{fig:summary}
\end{figure}

\cref{fig:summary} summarizes the unconstrained (in black) and constrained (in red) folds of the above two sub-sections. Each curve is the trajectory of $C'$ when we fix $a$ but vary $b$. The only difference between the black and red sets is the range of $b$: $[0,a]$ versus the more complex range implied by \cref{fig:region}.

\section{The rectangular page}
For a rectangular page, most of the previous results can be used. Given the page $[0,1]\times[0,A], A\ge1$, the only difference from the square page is the domain of $a$ increase from $[0,1]$ to $[0,A]$. 

For case 1, we again have 
$$
x_e(a,b)=x_{O'}(a,b) = \frac{2a^2b}{{a^2+b^2}}\quad 0\le b \le 1\le a.
$$
It is easy to show $\partial x_e/\partial b\ge 0$ and $\partial x_e/\partial a\ge 0$. So the maximum of $x_e$ is
$$
x_e(A,1)=\frac{2A^2}{{A^2+1}}.
$$
If $A>1$, $x_e>1$, so unlike the square page, this fold can be used as a bookmark (provided $A>1$). The optimal fold is clearly along the diagonal of the page. By the way, the maximum excess is
$$
e(A,1) = \frac{A^2-1}{{A^2+1}}.
$$
which means $\lim_{A\to\infty} e(A,1) = 1$, so this solution has a finite upper limit no matter now large $A$ is (i.e., how long the page is).

\begin{figure}[h]
	\centering
	\includegraphics[width=0.5\linewidth]{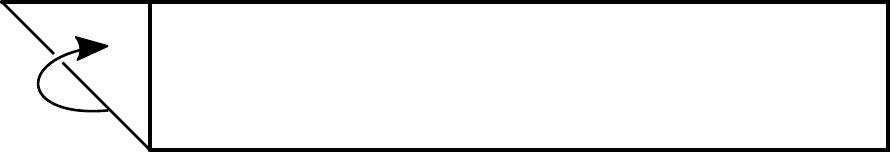}
	\caption{A solution that is close to optimal when the page is long and narrow (unconstrained problem). It is also the true optimal solution of the constrained problem when $A$ exceeds a critical value.}
	\label{fig:last}
\end{figure}

For case 2, we can reuse \cref{eqn:case2}, which is reproduced below:
\begin{equation}\label{eqn:case2rect}
e(a,b) = \frac{2b(a-b)}{1+(a-b)^2}.
\end{equation}
However, here the ranges of \emph{both} $a$ and $b$ are $[0,A]$. A contour plot of \cref{eqn:case2rect} indicates that the maximum of $e$ is still achieved at $a=A$, then setting $\dd e(A,b)/\dd b=0$ gives the optimal solution
\begin{equation}\label{eqn:case2rectsol}
e\left(A,\frac{1+A^2-\sqrt{1+A^2}}{A}\right) = \sqrt{1+A^2}-1
\end{equation}
It can be verified that setting $A=1$ recovers the optimal solution in \cref{sec:case2}. An interesting question is what happens when $A$ is very large. It is easy to show that \cref{eqn:case2rectsol} becomes approximately
\begin{equation}\label{eqn:case2rectsolapp}
e\left(A,A-1\right) = A-1
\end{equation}
This suggests that when the page is very long and narrow (e.g., the receipt of your grocery shopping), the solution is to fold such that the left edge overlaps with the top edge (see \cref{fig:last}). This solution does not exceed the top edge either, even though we never imposed this constraint!

\begin{figure}[h]
	\centering
	\includegraphics[width=0.5\linewidth]{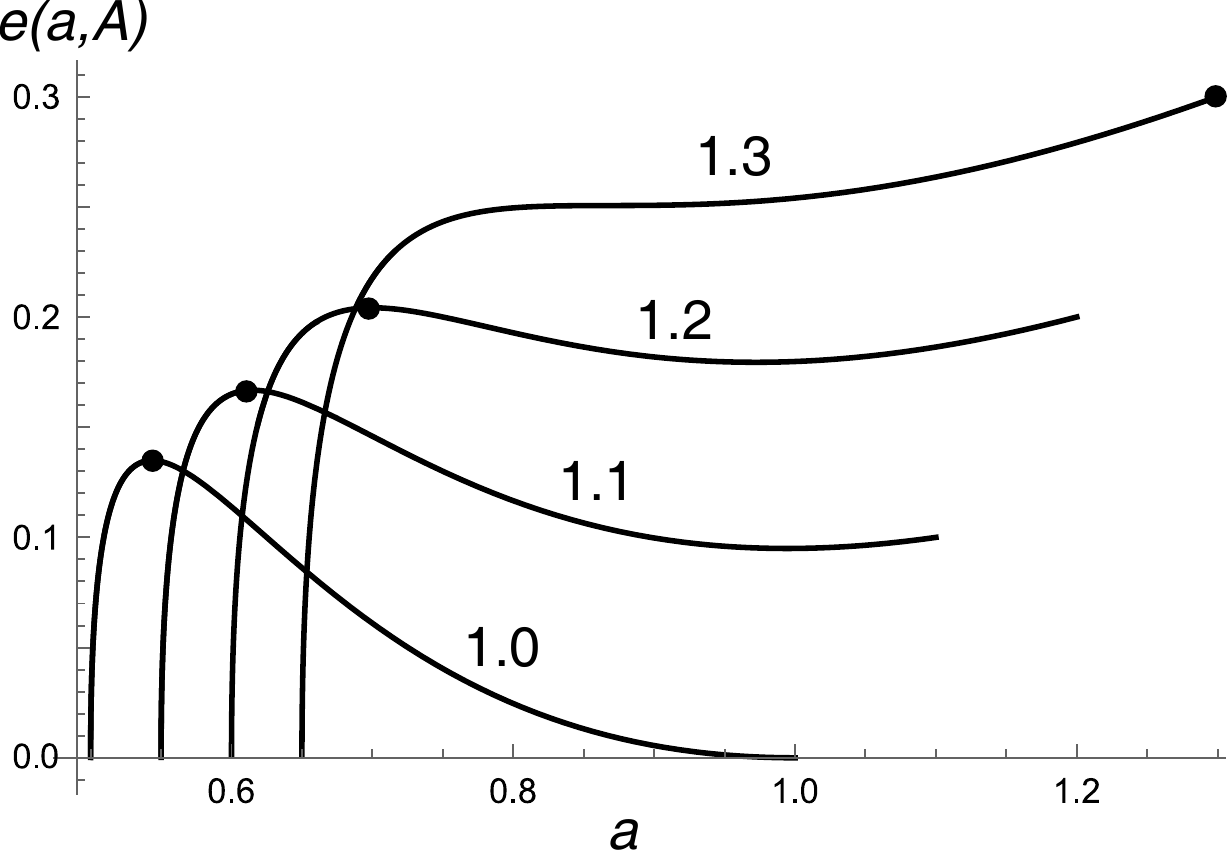}
	\caption{The curves $e(a,A)$ for four different $A$ values labeled by text. The black dots mark the $a$ values where the optimal $e$ is achieved.}
	\label{fig:transition}
\end{figure}

If we explicitly add $y_e\le A$ to the problem, the optimal solution still occurs at one point on the nonlinear boundary. However, here an interesting transition appears. When $A$ is relatively small, that (optimal) point is internal to the boundary; but when $A$ is larger, it coincides with vertex $(A,A-1)$ of the boundary. To shows this transition we can solve from the nonlinear boundary equation
$$\frac{2 a}{1+(a-b)^2}=A$$
and get (the other solution is discarded) 
$$
b=a-\frac{\sqrt{2 a A-A^2}}{A},\quad a\in\left[\frac{A}{2},A\right].
$$
Substitute into \cref{eqn:case2rect} gives the objective
$$
e(a,A) = \frac{A}{a}+\sqrt{A (2 a-A)}-2.
$$
Notice that the function is now of $a$ and $A$ instead of $a$ and $b$; the upper limit of $a$ or the page length has itself became a variable. \cref{fig:transition} shows four curves $e(a,A)$ with different $A$ values. For the first three curves, $A$ is relatively small, and $x_e$ occurs at $a<A$. For the last curve, the value of $e$ at the boundary point $a=A$ exceeds that of the internal point where $\dd e/\dd a=0$.

\begin{figure}[h]
	\centering
	\includegraphics[width=0.5\linewidth]{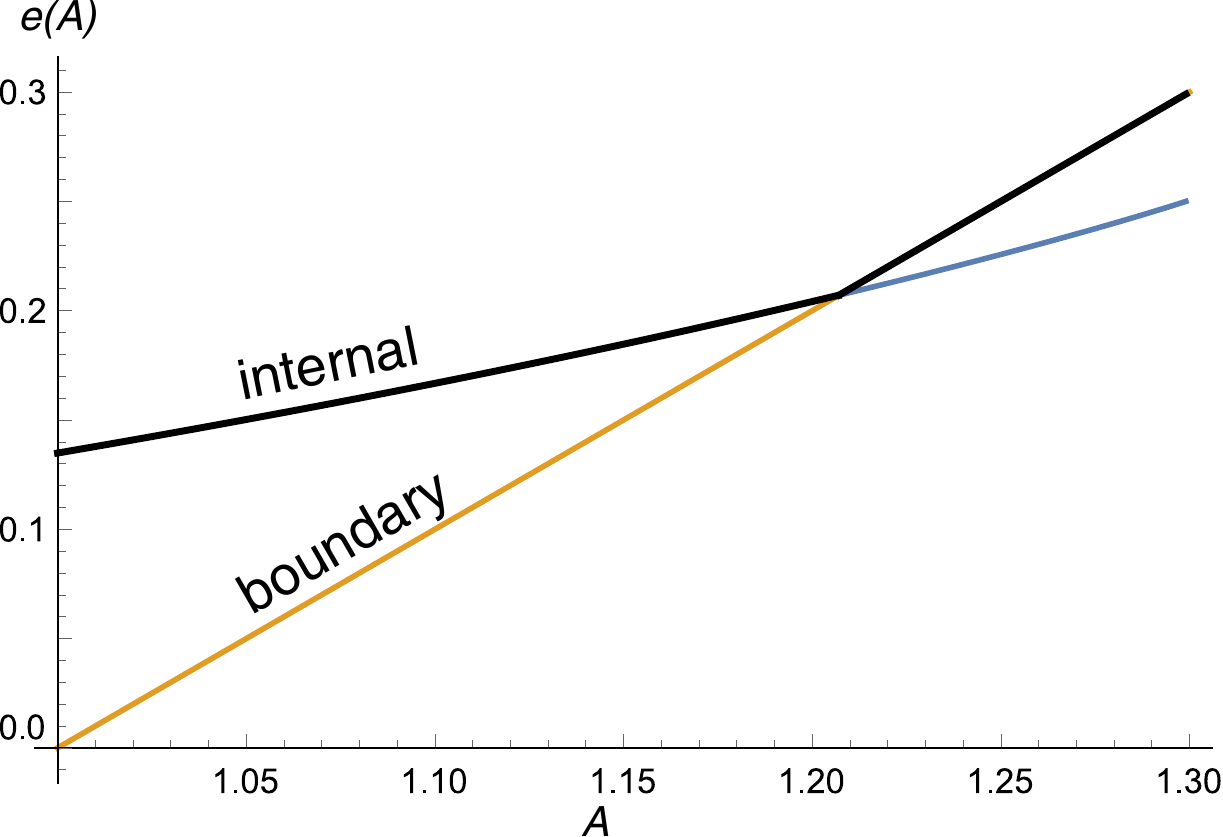}
	\caption{The ``phase diagram'' of $e(A)$ when the page is rectangular and the fold cannot exceed the top edge. The curve labeled by ``internal'' corresponds to the first three example curves in \cref{fig:transition}, whereas the one labeled by ``boundary'' corresponds to the fourth curve in that figure.}
	\label{fig:phase}
\end{figure}

The dependence of $e(A)$ on $A$ is more clearly shown by \cref{fig:phase}. When $A$ is small the optimal fold is ``irregular'' in that the crease is tilted at irrational angles (remember the solution in \cref{sec:foldconstraint}); when $A$ exceeds a certain critical value, the optimal fold becomes $e(A, A-1)$, which is highly ``regular'' (\cref{fig:last}). Interestingly, the second situation has the same solution as that of the unconstrained problem when $A$ is very large. 

The critical value of $A$ can be found numerically as $A_{cr}\approx1.20711$: passing it means a drastic change of the optimal folding method. The transition behavior of $e(A)$ is not unlike the ``first-order phase transition'' in thermodynamics, but here it is caused entirely by the geometric constraint. 

\section{Discussion}
This paper explores a simple problem: how to fold a page (square or rectangular) \emph{once} such that the folded part extends rightward to the maximum. We cast it into an optimization problem that can be either constrained or not depending on the needs. The key finding is that a rectangular page with constraint displays a nontrivial first-order phase transition behavior: the optimal way of folding it changes drastically when its aspect ratio exceeds a critical value of about 1.20711.

At least two directions of further studies can be thought of easily. The first direction departs from the context of notebook pages, and explores more complex shapes including circular or concave ones. Certainly, the extension to higher dimensions is also possible, although much harder to conceptualize. The second direction is to allow multiple folds: what if we can make two or more folds? The math is more complex, but at least we can easily show that the optimal solution of folding our square page with $y_e$ constraint can be beaten. This winning ``two-fold'' solution is as simple as first folding the square along diagonal $AC$ (\cref{fig:cases}, right panel), this generates an isosceles right triangle. The second step is to fold this triangle such that $AC$---now a triangle edge---aligns with the top edge. The excess $e=\sqrt{2}-1\approx0.414$ is greater than our optimal ``one-fold'' solution. However, whether it is optimal among all two-fold solutions and how it performs in a more general setting await further exploration.

\bibliographystyle{plainnat}
\bibliography{papers}

\begin{thebibliography}{5}
\providecommand{\natexlab}[1]{#1}
\providecommand{\url}[1]{\texttt{#1}}
\expandafter\ifx\csname urlstyle\endcsname\relax
  \providecommand{\doi}[1]{doi: #1}\else
  \providecommand{\doi}{doi: \begingroup \urlstyle{rm}\Url}\fi

\bibitem[Alperin and Lang(2009)]{alperin2009}
Roger~C. Alperin and Robert~J. Lang.
\newblock One-, two-, and multi-fold origami axioms.
\newblock \emph{Origami}, 4:\penalty0 371--393, 2009.

\bibitem[Hernandez et~al.(2018)Hernandez, Hartl, and Lagoudas]{hernandez2018}
Edwin A.~Peraza Hernandez, Darren~J. Hartl, and Dimitris~C. Lagoudas.
\newblock \emph{Active {{Origami}}: {{Modeling}}, {{Design}}, and
  {{Applications}}}.
\newblock {Springer}, July 2018.
\newblock ISBN 978-3-319-91866-2.

\bibitem[Hull(2012)]{hull2012}
Thomas Hull.
\newblock \emph{Project {{Origami}}: {{Activities}} for {{Exploring
  Mathematics}}, {{Second Edition}}}.
\newblock {CRC Press}, December 2012.
\newblock ISBN 978-1-4665-6809-9.

\bibitem[Lang(2003)]{lang2003}
Robert~J. Lang.
\newblock Origami and geometric constructions.
\newblock \emph{Self Published}, 2003.

\bibitem[Lang(2011)]{lang2011}
Robert~J. Lang.
\newblock \emph{Origami Design Secrets: Mathematical Methods for an Ancient
  Art}.
\newblock {AK Peters/CRC Press}, 2011.

\end{thebibliography}
\end{document}